\documentclass[11pt]{article}
\usepackage[latin1]{inputenc}
\usepackage[a4paper,scale={.82,.85}]{geometry}
\usepackage{amsmath,amssymb, amsfonts}
\usepackage{dsfont}

\usepackage[usenames,dvipsnames,svgnames,table]{xcolor}
\usepackage[colorlinks=true, pdfstartview=FitV,linkcolor=ForestGreen,citecolor=ForestGreen, urlcolor=blue]{hyperref}

\usepackage{theorem}%
{\theorembodyfont{\upshape}%
  \newtheorem{dfn}{Definition}
   \newtheorem{rmk}{Remark}
   \newtheorem{ex}{Example}
   \newtheorem{exs}{Examples}
}
{\theorembodyfont{\slshape}%
  \newtheorem{prop}{Proposition}
  \newtheorem{thm}[prop]{Theorem}
}
\newcommand{\cqfd}{\hfill\rule{0.35em}{0.35em}} 

\makeatother
\newcommand{\ie}{{\it i.e.} }		
\newcommand{\eg}[1][~]{{\it e.g.}#1}	

\newcommand{\N}{\mathds{N}}
  
\newcommand{\R}{\mathds{R}} 
\newcommand{\C}{\mathds{C}}
\newcommand{\HH}{\mathds{H}}

\newcommand{\norme}[2][]{\left\|#2\right\|_{#1}}
\let\div=\relax\DeclareMathOperator{\div}{div}

\begin{document}
\title{A simple proof of the Hardy inequality on Carnot groups\\
and for some hypoelliptic families of vector fields}
\author{F.~Vigneron\\[1ex]
\small\texttt{francois.vigneron@u-pec.fr}\\
\small Univ. Paris-Est Créteil, LAMA, UMR8050 du CNRS \\
\small 61, avenue du Général de Gaulle, F94010 Créteil, France.}
\maketitle
\begin{abstract}
We give an elementary proof of the classical Hardy inequality on any Carnot group, using only integration by parts and a fine analysis of the commutator structure,
which was not deemed possible until now.
We also discuss the conditions under which this technique can be generalized to deal with hypo-elliptic families of vector fields,
which, in this case, leads to an open problem regarding the symbol properties of the gauge norm.
\end{abstract}

The classical Hardy inequality \cite{HLP} on a smooth open domain $\Omega\subset \R^n$ ($n\geq 3$) reads:
\begin{equation}
\forall f\in H^1_0(\Omega), \qquad
\sup_{x_0\in\Omega}\left(
\int_{\Omega} \frac{|f(x)|^2}{|x-x_0|^2}\, dx\right)
\leq \frac{4}{(n-2)^2}\int_{\Omega}|\nabla f(x)|^2 dx.
\end{equation}
Since L.~D'Ambrosio~\cite{A1} it has been well known that similar inequalities hold on nilpotent groups, but interest on this
matter is still high; see \eg\cite{RS2017}, \cite{AM18}, \cite{APS2019}.

\medskip
An important reference on this matter is a recent note by~H.~Bahouri, C.~Fermanian, I.~Gallagher~\cite{BFG}.
It is dedicated to refined Hardy inequalities on graded Lie groups and relies on constructing a general Littlewood-Paley theory
and, as such,  involves the machinery of the Fourier transform on groups.

Expanding the generality towards hypoelliptic vector fields, G.~Grillo's~\cite{G} article contains an inequality that holds for $L^p$-norms,
without an underlying group structure, and contains weights that allow positive powers
of the Carnot-Caratheodory distance on the right-hand side.
However, the proof of this generalization involves the whole power of the sub-Riemannian
Calderon-Zygmund theory.

Another beautiful reference is the paper by P.~Ciatti, M.G.~Cowling, F.~Ricci~\cite{CCR2015} that studies these matters on stratified Lie groups,
but with the point of view of operator and interpolation theory (see also \cite{KMP99} to highlight some subtleties in this approach).

\medskip
The main goal of the present paper is to prove a general result on Carnot groups, albeit slightly simpler than those of \cite{BFG} or \cite{G},
by using only \textsl{elementary} techniques: most of the paper relies only on integrations by part and on a fine analysis of the commutator structure.
We will occasionally use interpolation techniques, but it is only required here if fractional regularities are sought after.

\medskip
A Carnot group is a connected, simply connected and nilpotent Lie group $G$ whose Lie algebra $\mathfrak{g}$ admits a stratification, \textit{i.e.}
\begin{equation}\label{DEFSTRAT}
\mathfrak{g} = \bigoplus_{j=1}^{m} V_{j}\qquad\text{where}\qquad [V_1,V_j]=V_{j+1}
\end{equation}
with $V_{m}\neq\{0\}$ but $[V_1,V_{m}]=\{0\}$.
The dimensions will be denoted by $q_j=\dim V_j$ and $q=\sum q_j = \dim\mathfrak{g}$.
Given a basis $(Y_\ell)_{\ell=1,\ldots,q}$ of $\mathfrak{g}$ adapted to the stratification each index~$i\in\{1,\ldots,q\}$ can be associated to a unique \textsl{weight}
$\omega_\ell\in\{1,\ldots,m\}$ such that $Y_\ell\in V_{\omega_\ell}$, namely
\begin{equation}\label{def:njomegaj}
\omega_\ell = j \qquad\text{for}\qquad n_{j-1}<\ell\leq n_j
\end{equation}
where $n_0=0$ and $n_{j}=n_{j-1}+q_j$ for $j=1,\ldots,m$ is the sequence of cumulative dimensions. Note that $n_1=q_1$ and $n_m=q$.
 The horizontal derivatives are the derivatives in the first layer (see \S\ref{S:horDer} below) and they are collected together in the following notation:
\begin{equation}\label{nablaG}
\nabla_G f=(Y_1^Lf,\ldots,Y_{n_1}^Lf).
\end{equation}
The stratification hypothesis ensures that each derivative $Y_if$ can be expressed as at most $\omega_i-1$ commutators of horizontal derivatives.
The \textsl{homogeneous dimension} is the integer
\begin{equation}\label{HOMDIM}
Q=\sum_{j=1}^m j q_j=\sum_{\ell=1}^q\omega_\ell.
\end{equation}
For $k\in\N$, the Sobolev space $H^k(G)$ is the subspace of functions $\phi\in L^2(G)$  such that $\nabla_G^\alpha \phi \in L^2(G)$ for any multi-index $\alpha$
of length $|\alpha|\leq k$. Fractional spaces can, for example, be defined by interpolation.
The main result that we intend to prove here is the following.

\begin{thm}\label{main_thm}
Let $G$ be a Carnot group and $\norme{\cdot}$ any homogeneous pseudo-norm equivalent to the Carnot-Caratheodory distance to the origin. Then, for any
real $s$ with $0\leq s< Q/2$, there exists a constant $C_s>0$ such that:
\begin{equation}\label{MAINHARDY}
\int_G \frac{|f(g)|^2}{\norme[G]{g}^{2s}} dg \leq C_s \norme[H^s(G)]{f}^2
\end{equation}
for any function $f\in H^s(G)$.
\end{thm}
A similar Hardy inequality was proved by the author in~\cite{PhD} for families of vector fields that satisfy
a H\"ormander bracket condition of step 2; the proof was based on the ideas of \cite{BCG}  and \cite{BahCoh} but was never published independently.
This result was part of a broader study \cite{BCX1}, \cite{BCX2}, \cite{MV07}, \cite{V07} aiming at characterizing the traces of Sobolev spaces on the Heisenberg group, along
hypersurfaces with non-degenerate characteristic points.
Here, instead, we concentrate (except in~\S\ref{S:vectFields}) on the case of stratified groups, but without restrictions on the step $m$ of the stratification.

\medskip
The mathematical literature already contains numerous Hardy-type inequalities either on the Heisenberg group, for
the $p$-sub-Laplacian, for Grushin-type operators and $H$-type groups  (see \eg\cite{A1}, \cite{Kombe}).
Sometimes (\eg in \cite{GL}, \cite{N1}, \cite{N4}), a weight is introduced in the left-hand side that
vanishes along the center of the group, \textit{i.e.} along the (most) sub-elliptic direction.
For example,~\cite{GL} contains the following inequality
on the Heisenberg group~$\HH_n\simeq \C^n\times\R$:
\begin{equation}
\int_{\HH_n} \frac{|f(x)|^2}{d_{\HH_n}(x,0)^2} \:\Phi(x)dx \leq A\sum_{j=1}^n
\left(\norme[L^2]{X_jf}^2 +\norme[L^2]{Y_j f}^2\right) +B \norme[L^2]{f}^2
\end{equation}
where $(X_j,Y_j)$ are a basis of the first layer of the stratification and $d_{\HH_n}((z,t),0)\simeq \sqrt[4]{|z|^4+t^2}$ is the gauge
distance and $\Phi$ is a cut-off function that vanishes along the center $z=0$:
$$\Phi(z,t)=\frac{|z|^2}{\sqrt{|z|^4+t^2}}\cdotp$$
A secondary goal of this article is to show that such a cut-off is usually not necessary.

\medskip
The core of our proof of theorem~\ref{main_thm} (see \S\ref{PROOF}) consists in an integration by part against the radial field (the infinitesimal generator of dilations).
The radial field can be expressed in terms of the left-invariant vector fields but, as all strata are involved, this first step puts $m-s$ too many derivatives
on the function. Next, one uses the commutator structure of the left-invariant fields to carefully backtrack all but one derivative and let them act instead on the coefficients
of the radial field. This step requires that those coefficients have symbol-like properties. One can then conclude by an iterative process that reduces the
Hardy inequality with weight $\norme[G]{g}^{-s}$ to the one with weight $\norme[G]{g}^{-(s-1)}$ as long as $s<Q/2$.

\medskip
Finally, it is worth mentioning that a byproduct of our elementary approach concerns the symbol properties of the Carnot-Caratheodory norm
(or of any equivalent gauge). For general hypoelliptic families of  vector fields, \textbf{the norm is not always a symbol of order~1} (see section~\S\ref{S:vectFields}).
On the contrary, on Carnot groups, it happens to  always be equivalent to such a symbol (Proposition~\ref{prop_symbol}).
At the end of the article, we discuss sufficient conditions for this property to hold for families of hypoelliptic vector fields, based either on
the order $m$ of the H\"ormander condition (Theorem~\ref{STEP2FIELDS}), or on the way the commutators are structured (Theorems~\ref{ABSTRACTSTEPmFIELDS} and~\ref{STEP3FIELDS}).

\medskip
The structure of the article goes as follows. Section \ref{S:Intro} is a brief survey of calculus on
Carnot groups. It also sets the notations used subsequently. Section~\ref{S:proof} contains the actual proof of theorem~\ref{main_thm}
and concludes on  theorem~\ref{hom_thm}, which is the homogeneous variant of the previous statement.
Section~\ref{S:vectFields} adresses an open question regarding families of vector fields that satisfy a H\"ormander
bracket condition, but lack an underlying group structure.

\medskip
I would like to express my gratitude to J.-Y.~Chemin, who brought this problem to my attention a long time ago,
and to wish him a happy 60th birthday.

\section{A brief survey of calculus on Carnot Groups}\label{S:Intro}

Let us first recall some classic definitions and facts about nilpotent Lie groups.
We also introduce notations that will be needed in \S\ref{S:proof}.
For a more in-depth coverage of Lie groups, sub-Riemannian geometry and nilpotent groups,
see~\eg\cite{M}, \cite{FS}, \cite{R1} or the introduction of~\cite{rigot}.

\subsection{Left-invariant vector fields and the exponential map}

Let us consider a Lie group $G$ and $\mathfrak{g}=T_eG$ its Lie algebra;
$e$ denotes the unit element of~$G$. Left-translation is defined by $L_g(h)=gh$.

\begin{dfn} A vector field $\xi$ is called \textsl{left-invariant}
if $(L_g)_\ast \circ \xi = \xi\circ L_g$. 
Such a vector field is entirely determined by $v=\xi(e)\in\mathfrak{g}$. To signify that $v$ generates $\xi$,
one writes $\xi=v^L$ thus:
\begin{equation}\label{LEFTINVFIELD}
v^L(g)=d(L_g)_{\vert e} (v).
\end{equation}\end{dfn}
The tangent bundle $TG$ identifies to $G\times\mathfrak{g}$ by the map $(g,v) \mapsto (g,v^L(g))$.

\bigskip
The flow $\varPhi_t^v$ of a left-invariant vector field $v^L$ exists for all time. Indeed,
one has $\varPhi_t^v(g) = L_g\circ \varPhi_t^v(e)$, which implies that
$\varPhi_{t+s}^v(e)=L_{\varPhi_s^v(e)}\circ \varPhi_t^v(e)$,
thus allowing the flow to be extended globally once it has been constructed locally.
\begin{dfn}
The \textsl{exponential map} $\exp:\mathfrak{g}\to G$ is defined by $\exp(v)=\varPhi_1^v(e)$
where $(\varPhi_t^v)_{t\in\R}$ is the flow of the left-invariant vector field $v^L$.
\end{dfn}
One can check that the flow of $v^L$ starting from $g\in G$ is $\varPhi_t^v(g) = g \exp(tv)$.
In particular, $$\exp(sw)\exp(tv) = \Phi_t^v(\exp(sw)) \qquad\text{and}\qquad
d\exp_{\vert 0}=\operatorname{Id}_\mathfrak{g}.$$

\subsection{The Baker-Campbell-Hausdorff formula}

The commutator of two left-invariant vector fields is also a left-invariant field. Therefore,
the commutator of $u,v\in\mathfrak{g}$ is defined by $[u,v]=[u^L,v^L](e)\in\mathfrak{g}$.
The product law of $G$ induces an extremely rigid relation between exponentials, known as
the Baker-Campbell-Hausdorff formula:
\begin{equation}\label{BCH}
\exp(u)\exp(v)=\exp(\mu(u,v))
\end{equation}
where $\mu(u,v)=u+v+\frac{1}{2}[u,v]+\frac{1}{12}\left([u,[u,v]]+[v,[v,u]]\right)+\ldots{}$ is a
universal Lie serie in $u, v$ \ie an expression consisting of the iterated commutators of $u$ and $v$.
In general,
this formula holds provided $u$ and $v$ are small enough for the series to converge (see \eg\cite[\S1.3]{R1}).
Subsequently, one will only use the linear part of \eqref{BCH} with
respect to one variable:
\begin{equation}\label{BCH:LIN}
d\mu(u,\cdot)_{\vert 0}(w)
= \frac{\operatorname{ad}(u)}{1-e^{-\operatorname{ad}(u)}}(w)
=w+\sum_{n=1}^\infty  \frac{(-1)^nB_n}{n!} \:
\underset{n\text{ times }u}{[u,\ldots,[u},w]]
\end{equation}
where $B_n$ are the Bernoulli numbers (\ie$\frac{x}{e^x-1}=\sum \frac{B_n}{n!}x^n$)
and $\operatorname{ad}(u)=[u,\cdot]$.
This formula is classical and can be found \eg in \cite{R1} or \cite{KO}.

\subsection{Stratification}

From now on, $G$ is supposed to be \textsl{stratified}, \ie it is a Carnot group as defined in the introduction of this paper.
A stratified group is, in particular, nilpotent of step $m$. Moreover, elementary linear algebra gives restrictions on the possible dimensions
$q_j=\dim V_j$ of the strata:
$$q_2 \leq \frac{q_1(q_1-1)}{2} \qquad\text{and for }j\geq2,\qquad q_{j+1}< q_1q_j.$$
The last inequality is strict because of the Jacobi identity $[u,[v,w]]-[v,[u,w]]=-[w,[u,v]]$.
For an exact count of the possible relations, 
see \eg\cite{R2}.

\begin{prop}\label{P:globDiff}
Let $G$ be a Carnot group. Then $\exp:\mathfrak{g}\to G$ is a global diffeomorphism that allows $G$
to be identified with the set $\mathfrak{g}$ equipped with the
group law $u\ast v=\mu(u,v)$. The identity element  is $0$ and the inverse of $u$ is $-u$.
\end{prop}

\paragraph{Proof.}
This claim is very standard so one only sketches the proof briefly.
As $d\exp_{\vert 0}=\operatorname{Id}_\mathfrak{g}$, there exists a neighborhood $U_0$ of $e$ in $G$ and
$V_0$ of $0$ in $\mathfrak{g}$ such that $\exp : V_0 \to U_0$ is a diffeomorphism.  
As $G$ is connected, it is generated by any neighborhood of $e$ and in particular by $U_0=\exp(V_0)$.
But, as $\mathfrak{g}$ is nilpotent, the expression $\mu(u,v)$ is a Lie polynomial of order~$m$,
thus \eqref{BCH} holds for any $u,v\in \mathfrak{g}$. Combining these facts implies
that the exponential map is surjective.
Next, one can show that the pair $(\mathfrak{g},\exp)$ is a covering space of $G$.
Indeed, given $g=\exp(v)\in G$,
one gets a commutative diagram of diffeomorphisms:
$$\begin{array}{rcl}
V_0 & \overset{\exp}{\longrightarrow} & U_0\\
\mu(v',\cdot)\Big\downarrow & & \Big\downarrow L_{g}\\
V' & \underset{\exp}{\longrightarrow} & g\cdot U_0
\end{array}$$
for each $v'\in\mathfrak{g}$ such that $\exp(v')=g$.
Finally, by a standard covering space argument based on the fact that $\mathfrak{g}$ is path connected (as vector space) and $G$ is simply connected (in the stratification assumption), one can claim that
the exponential map is a global diffeomorphism.
\cqfd

\begin{ex}
The Heisenberg group $\HH$ can be realized as a set of upper-triangular matrices with diagonal entries equal to~$1$.
The group law in $\HH$ is the multiplication of matrices:
\[
\begin{pmatrix} 1 & p & r\\0 & 1 & q\\0&0&1\end{pmatrix}
\cdot
\begin{pmatrix} 1 & p' & r'\\0 & 1 & q'\\0&0&1\end{pmatrix}
=
\begin{pmatrix} 1 & p+p' & r+r'+pq'\\0 & 1 & q+q'\\0&0&1\end{pmatrix}.
\]
The Lie algebra of $\HH$ is
\[
\mathfrak{h}=\left\{\begin{pmatrix} 0 & p & r\\0 & 0 & q\\0&0&0\end{pmatrix}=pY_1+qY_2+rY_3 \,;\, p,q,r\in\R\right\}.
\]
Left-invariant vector fields on $\HH$ are linear combinations of
\[
Y_1^L(g)=Y_1, \quad Y_2^L(g)=Y_2+pY_3 \quad\text{and}\quad Y_3^L(g)=Y_3
\quad\text{where}\quad
g=\begin{pmatrix} 1 & p & r\\0 & 1 & q\\0&0&1\end{pmatrix}.
\]
$\HH$ is a stratified nilpotent group with $V_1=\operatorname{Span}(Y_1,Y_2)$ and $V_2=\operatorname{Span}(Y_3)$.
The exponential map is the usual exponential of nilpotent matrices. It transfers the group structure to $\mathfrak{h}\simeq\R^3$ by~\eqref{BCH}:
\[
(p,q,r)_\mathfrak{h} \simeq \begin{pmatrix} 0 & p &r\\0 & 0 & q\\0&0&0\end{pmatrix},
\qquad
\exp \begin{pmatrix} 0 & p & r\\0 & 0 & q\\0&0&0\end{pmatrix}
=\begin{pmatrix} 1 & p & r+\frac{1}{2}pq\\0 & 1 & q\\0&0&1\end{pmatrix},
\]
\[
\mu\left( (p,q,r)_\mathfrak{h},(p',q',r')_\mathfrak{h} \right)=\left(p+p',q+q',r+r'+{\textstyle \frac{1}{2}}(pq'-qp')\right)_\mathfrak{h}.
\]
In exponential coordinates on $\mathfrak{h}\simeq \R^3$, the left-invariant vector fields thus take the following form:
\[
\tilde Y_1^L(p,q,r) = \partial_p -\frac{1}{2} q\partial_r, \qquad
\tilde Y_2^L(p,q,r) = \partial_q +\frac{1}{2} p\partial_r, \qquad
\tilde Y_3^L(p,q,r) = \partial_r.
\]
The two different expressions of the fields correspond to the change of variables $(p,q,r)\mapsto (p,q,r-\frac{1}{2}pq)$.
\end{ex}

\paragraph{Remarks}
\begin{itemize}
\item
In general (even if $G$ is connected and nilpotent), only the Lie group action can be recovered from the exponential map but not $G$ itself.
For example, $G=\{z\in\C\,;\, |z|=1\}\simeq \mathbb{S}^1$ with rotation law $(z,z')\mapsto zz'$ is a nilpotent group.
One has $\mathfrak{g}=\R$ and $\mu(x,y)=x+y$ but the exponential map is $\exp(x)=e^{ix}$ and is obviously not
a global diffeomorphism.
\item 
Combined with~\eqref{BCH:LIN}, the commutative diagram of the proof of Proposition~\ref{P:globDiff}
provides a general formula for the differential of the exponential map, which we will need later on.
For $v,w\in\mathfrak{g}$ and $g=\exp(v)$, one has:
\begin{equation}\label{DIFFEXP}
d \exp_{\vert v} (w) = 
 \sum_{n=0}^{m-1} \frac{(-1)^n}{(n+1)!}\:\underset{n\text{ times }v}{[v,\ldots,[v},w]]^L(g).
\end{equation}
For example, on the Heisenberg group, one gets $d \exp_{\vert v} (w)=(w-\frac{1}{2}[v,w])^L(g)$.
\end{itemize}
\par

\begin{prop}
For any indices $j,k\in\{1,\ldots,m\}$, one has
\begin{equation}\label{HIGHCOMM}
[V_j,V_k]\subset \begin{cases}
V_{j+k} & \text{if }j+k\leq m,\\
\{0\} &\text{otherwise.}
\end{cases}
\end{equation}
\end{prop}
\paragraph{Proof.}
By convention, let us write $V_n=\{0\}$ if $n>m$.
For $k=1$, the property holds by definition. For $k=2$, as $V_2=[V_1,V_1]$, any element can be written $[X,Y]$ with $X,Y\in V_1$. For $Z\in V_j$, one uses the identity $[A,BC]=[A,B]C+B[A,C]$ to get:
\[
[Z,[X,Y]]=[[Z,X],Y]-[[Z,Y],X]\in [[V_j,V_1],V_1]\subset [V_{j+1},V_1]\subset V_{j+2}.
\]
Next, one proceeds recursively. Assuming that for some $k\geq2$, one has $[V_j,V_k]\subset V_{j+k}$ for any $j$, then given $Z\in V_{k+1}=[V_1,V_k]$, one writes $Z=[X,\zeta]$ with
$X\in V_1$ and $\zeta\in V_k$. Then for any $W\in V_j$, the Jacobi identity gives
\[
[W,Z]=[W,[X,\zeta]]=-[X,[\zeta,W]]-[\zeta,[W,X]]\in [V_1,V_{j+k}]+[V_k,V_{j+1}] \subset V_{j+k+1}
\]
which makes the property hereditary in $k$.
\cqfd

\subsection{Stratified dilations}

The next essential object in a Carnot group is the \textsl{dilation} of the Lie algebra:
\begin{equation}\label{AlgebraDILATION}
\forall r>0,\qquad
\delta_r = \sum_{j=1}^m r^j \pi_j
\end{equation}
where $\pi_j:\mathfrak{g}\to V_j$ is the projection onto $V_j$ with kernel $\bigoplus_{k\neq j} V_k$.
Identifying $G$ to $\exp(\mathfrak{g})$, one gets a one parameter family of group automorphisms that
we will simply denote by 
\begin{equation}\label{DILATION}
r g=\exp\circ \,\delta_r \circ \exp^{-1}(g)
\end{equation}
for any $r>0$ and $g\in G$.

\bigskip
The next result is an immediate consequence of the definition but should later be compared with the scaling property~\eqref{RSCALING} of the radial vector field.
\begin{prop}
The dilation of a left-invariant vector field $v^L$ is given by:
\begin{equation}
(\delta_r v)^L(r g)=\left((L_{r g})_\ast\circ\,\delta_r\circ(L_g)_\ast^{-1}\right)\left(v^L(g)\right).
\end{equation}
\end{prop}

\bigskip
Up to a constant factor,
the Haar measure on~$G$ is given by the Lebesgue measure on $\mathfrak{g}\simeq \R^q$ and is commonly denoted $dg$.
\begin{prop}
One has:
\begin{equation}
\forall\phi\in L^1(G),\quad\forall r>0, \qquad  \int_G \phi(r g) dg = r^{-Q}\int_G \phi(g)dg
\end{equation}
where $Q$ is the homogeneous dimension \eqref{HOMDIM} of $G$.
Note that when $m\neq1$, one has $Q>q$.
\end{prop}

\bigskip
Given (arbitrary\footnote{For a given basis of $\mathfrak{g}$ adapted to the stratification, one will chose here a Euclidian structure that renders this basis orthonormal. This is the natural choice when one proceeds to the identification $\mathfrak{g}\simeq\R^q$ through this basis.}) Euclidian norms $\norme[V_j]{\cdot}$ on each $V_j$ and $w=2\operatorname{LCM}(1,\ldots,m)$,
the anisotropic \textsl{gauge-norm} of either $v\in\mathfrak{g}$ or of $g=\exp(v)\in G$ is defined by:
\begin{equation}\label{GAUGE0}
\norme[\mathfrak{g}]{v} = \norme[G]{g} = \left(\sum_{j=1}^m \norme[V_j]{\pi_j(v)}^{w/j}\right)^{1/w}.
\end{equation}
The gauge norm is homogeneous in the following sense:
\begin{equation}
\norme[G]{r g}=r\norme[G]{g} \qquad \text{and}\qquad |\Box_{g,r}|=c_0 r^Q
\end{equation}
with a uniform constant $c_0$ and where the gauge-ball is defined by $$\Box_{g,r}=\{h\in G\,;\, \norme[G]{h^{-1}g}< r\}.$$

\begin{rmk} The intrinsic metric objects of $G$ are the so-called Carnot balls defined as the set
of points that can be connected to a center~$g_0\in G$ by an absolutely continuous path $\gamma$ whose velocity
is sub-unitary for almost every time \ie~such that $\dot\gamma(t)\in  (L_{\gamma(t)})_\ast(B_0)$
where $B_0\subset V_1$ is a fixed Euclidian ball of the first layer of the stratification (up to some
choice of a Euclidian metric on $V_1$). However, the ball-box theorem \cite{M} states that such intrinsic objects
can be sandwiched between two gauge-balls of comparable radii. For our purpose (the analysis of Sobolev spaces), one can  thus
deal only with gauge-balls without impeding  the generality. 
\end{rmk}

\subsection{Horizontal derivatives and Sobolev spaces on $G$}\label{S:horDer}

Vector fields $\xi$ on $G$ are identified with derivation operators on $C^\infty(G)$ by the Lie derivative
formula:
\begin{equation}
(\xi \phi)(g) = d\phi_{\vert g}(\xi(g)).
\end{equation}

\begin{dfn}
The \textsl{horizontal derivatives} are the left-invariant vector-fields associated with $V_1$.
\end{dfn}
Let us consider a basis $(Y_\ell)_{1\leq \ell\leq q}$ of $\mathfrak{g}$ that is adapted to the stratification, \ie:
\begin{equation}
V_j = \operatorname{Span}\: (Y_\ell)_{n_{j-1}< \ell\leq n_j}
\end{equation}
where $n_j$ is defined by~\eqref{def:njomegaj}.
An horizontal derivative is thus a vector field
$$\xi =\sum_{j\leq n_1} \alpha_j Y_j^L=\alpha\cdot\nabla_G$$
where $\alpha=(\alpha_1,\ldots,\alpha_{n_1})\in \R^{n_1}$ and
$\nabla_{G}$ is defined by~\eqref{nablaG}.
Non-commutative multi-indices are defined as follows:
for $\gamma=(\gamma_1,\ldots,\gamma_\ell)\in\{1,\ldots,n_1\}^\ell$, one writes $l=|\gamma|$ and $\nabla_G^\gamma=Y_{\gamma_1}^L\circ\cdots\circ Y_{\gamma_l}^L$.

\bigskip\noindent
Let us unfold the commutator structure in $\mathfrak{g}$ with the following notation:
\begin{equation}\label{COMMSTRUCT}
[Y_{\ell_1},\ldots,[Y_{\ell_n},Y_{\ell_{n+1}}]]=\sum_{\ell'} \kappa(\ell_1,\ldots,\ell_n,\ell_{n+1};\ell') Y_{\ell'}.
\end{equation}
Note that according to \eqref{HIGHCOMM}, one can warrant that $\kappa(\ell_1,\ldots,\ell_n,\ell_{n+1};\ell')=0$
if $\omega_{\ell'}\neq \omega_{\ell_1}+\ldots+\omega_{\ell_{n+1}}$.

\begin{rmk} To simplify computations, one can always assume that the basis is chosen such that:
\begin{equation}\label{COMMREVSTRUCT}
\forall \ell\in\{1,\ldots,q\}, \qquad 
Y_\ell = [Y_{\alpha_1(\ell)},\ldots,[Y_{\alpha_{k-1}(\ell)},Y_{\alpha_k(\ell)}]]
\end{equation}
where $k=\omega_\ell$ and $\alpha_i(\ell)\leq n_1$. Indeed, the Lie algebra is linearly generated by the commutators of the restricted family $\nabla_G$ and one just has to extract a basis from it.
\end{rmk}

\begin{ex} On the Heisenberg group $\HH$, the horizontal derivatives are 
left-invariant vector fields of the form $\xi=(\alpha Y_1 + \beta Y_2)^L$
for $\alpha,\beta\in\R$.
\end{ex}

\begin{dfn}
For $s\in\N$, the Sobolev space $H^s(G)$ consists of the functions such that each composition of at most~$s$ horizontal derivatives belongs to $L^2(G)$.
The norm is defined (up to the choice of the $Y_\ell$) by:
\begin{equation}\label{norme}
\norme[H^s(G)]{\phi}^2 = \sum_{|\gamma|\leq s} \int_G |\nabla_G^\gamma\phi(g)|^2 dg
\end{equation}
with $\gamma$ a non-commutative multi-index.
\end{dfn}

\begin{rmk}
The space $H^{2}(G)$ is the domain of the hypo-elliptic Laplace operator 
\begin{equation}\label{HYPOLAPLACE}
\mathcal{L}_G=-\sum_{\ell\leq n_1} (Y_\ell^L)^\ast\: Y_\ell^L.
\end{equation}
A celebrated result of L.~H\"ormander \cite{hormander} states that $H^s(G)\subset H^{s/m}_{\text{loc}}(\R^q)$
where the last Sobolev space is the classical one (homogeneous and isotropic) on $\R^q$.
\end{rmk}

\subsection{Exponential coordinates on a stratified group}\label{S:expCoord}

Given a basis $(Y_\ell)_{1\leq \ell\leq q}$ of $\mathfrak{g}$ adapted to the stratification,
one can define a natural coordinate system on $G$, called
\textsl{exponential coordinates}. Given $g=\exp(v)\in G$, its coordinates $x(g)=(x_\ell(g))_{1\leq\ell\leq q}\in\R^q$ are defined by:
\begin{equation}
v=\sum_{\ell=1}^q x_\ell(g) Y_\ell.
\end{equation}
The projections $(\pi_j)_{1\leq j\leq m}$ introduced in \eqref{AlgebraDILATION} are:
\begin{equation}\label{PROJECTIONS}
\forall j\in\{1,\ldots,m\},\qquad
\pi_j(v) = \sum_{\ell=1+n_{j-1}}^{n_j} x_\ell(g)Y_\ell.
\end{equation}
In exponential coordinates, the expression of stratified dilations \eqref{DILATION} is:
\begin{equation}
\forall \ell\in\{1,\ldots,q\},\qquad
x_\ell(r g) = r^{\omega_\ell} x_\ell(g).
\end{equation}
The gauge norm \eqref{GAUGE0} is given (for some fixed large $w\in\N$) by:
\begin{equation}\label{GAUGE1}
\norme[G]{g}=\norme[\mathfrak{g}]{x(g)}=\left(\sum_{j=1}^m \left(\sum_{\ell=1+n_{j-1}}^{n_j}|x_\ell(g)|^2\right)^{w/j} \right)^{1/(2w)}.
\end{equation}
One could however take any uniformly equivalent quantity as a gauge norm, which will be the case subsequently,
after proposition~\ref{prop_symbol}.

\begin{ex}
With the previous notations, the exponential coordinates on the Heisenberg group $\HH$ are:
\[
x_1(g)=p, \quad x_2(g)=q, \quad x_3(g) = r-\frac{1}{2}pq
\quad\text{for}\quad
g=\begin{pmatrix} 1 & p & r\\0 & 1 & q\\0&0&1\end{pmatrix}\in\HH.
\]
\end{ex}

\subsection{Left-invariant basis of vector fields}

When doing explicit computations, it is natural to identify $\mathfrak{g}$ with $\R^q$ through the previous coordinates. Given $v\in \mathfrak{g}$,
the left-invariant vector field $v^L$ on $G$ defined by~\eqref{LEFTINVFIELD} matches a corresponding vector field on $\mathfrak{g}\simeq\R^q$
that we will denote  by $\tilde{v}^L$.
According to the Baker-Campbell-Hausdorff formula~\eqref{BCH:LIN}:
\begin{equation}\label{BCH2}
\forall x=(x_1,\ldots,x_q)\in\mathfrak{g},\qquad
\tilde{v}^L(x)=v + \sum_{n=1}^{m-1}\sum_{\ell_1,\ldots,\ell_n} \frac{(-1)^n B_n}{n!}x_{\ell_1}\ldots x_{\ell_n} [Y_{\ell_1},\ldots,[Y_{\ell_n},v]]
\end{equation}
where each of the $\ell_i$ ranges over $\{1,\ldots,q\}$.

\medskip
After the identification $\mathfrak{g}\simeq\R^q$ and to avoid confusion, let us denote $Y_\ell$ by $\partial_\ell$:
the vectors $(\partial_\ell)_{1\leq \ell\leq q}$ are the dual basis of the coordinates $(x_\ell)_{1\leq\ell\leq q}$.
The left-invariant basis then becomes explicit:
\begin{equation}\label{HDERIVATIVE}
\forall x=(x_1,\ldots,x_q)\in\mathfrak{g},\qquad
\tilde{Y}_\ell^L(x)=\partial_\ell + \sum_{\ell'=1}^q \zeta_{\ell,\ell'}(x_1,\ldots,x_q)\partial_{\ell'}
\end{equation}
with, thanks to~\eqref{BCH2} and~\eqref{COMMSTRUCT}:
\begin{equation}\label{ZETA}
\zeta_{\ell,\ell'}(x) = \sum_{n=1}^{m-1} \frac{(-1)^n B_n}{n!} \:
\sum_{\ell_1,\ldots,\ell_n} \kappa(\ell_1,\ldots,\ell_n,\ell;\ell') \, x_{\ell_1}\ldots x_{\ell_n}.
\end{equation}
Let us point out that $\zeta_{\ell,\ell'}=0$ if $\omega_{\ell'}\leq \omega_\ell$ (because $\kappa$ vanishes),
thus the left-invariant correction to~$\partial_\ell$ only involves derivatives of a strictly higher weight.
In other words, the indices in~\eqref{HDERIVATIVE} can be restricted to $\ell'>n_{\omega_\ell}=q_1+q_2+\ldots+q_{\omega_\ell}$.

\medskip
Note also that $\zeta\in C^\infty(\R^q,\mathcal{M}_{q,q}(\R))$ and $\zeta(0)=0$.
More precisely, this matrix represents the differential action of left-translations, expressed in exponential coordinates:
\begin{equation}
(dL_g)_{\vert e} \equiv \operatorname{Id}_{\R^q} +\,\zeta(x(g)).
\end{equation}
Moreover, as $|x_\ell(g)|\lesssim \norme[G]{g}^{\omega_\ell}$, one has $|\zeta_{\ell,\ell'}(x(g))|\lesssim \norme[G]{g}^{\omega_{\ell'}-\omega_\ell}$.

\section{Proof of Theorem~\ref{main_thm}}\label{S:proof}

This section is devoted to the proof of the main statement. The key idea is to prove the
result for $s=1$ and then ``push'' the result up to the maximal regularity using only integrations by part.
Adding an interpolation step once the result is known for $s=1$, but before pushing it to a higher regularity,
allows one to capture all eligible fractional derivatives.
The actual proof is written in the last subsection~\S\ref{PROOF} but some preliminary results are required.

\subsection{Symbol classes $S^\alpha_n(G)$}

Symbol classes are a convenient way to classify the coefficients involved in the computations in terms of how they vanish at the origin.
\begin{dfn}
For $\alpha\in\R_+$ and $n\in\N\cup\{\infty\}$, the symbol class $S^\alpha_n(G)$ is defined as the set of functions $\phi\in L^\infty_{\text{loc}}(G)$ 
such that for any multi-index $\gamma$ of length $|\gamma|\leq n$, there exists
a constant $C_\gamma>0$ that ensures the following inequality:
\begin{equation}
\forall g\in G, \qquad \norme[G]{g}\leq 1 \qquad\Longrightarrow\qquad
|\nabla_G^\gamma \phi(g)| \leq C_\gamma \norme[G]{g}^{(\alpha-|\gamma|)_+}.
\end{equation}
For example, the symbols of class $S^0_0(G)=L^\infty_{\text{loc}}(G)$ are only required to be bounded near the origin.
The symbol class $S^\alpha_\infty(G)$ is also denoted $S^\alpha(G)$.
\end{dfn}
The following properties hold.
\begin{enumerate}
\item
The Leibnitz formula gives
$$\phi\in S^\alpha_m(G), \quad \psi\in S^\beta_n(G) \qquad\Longrightarrow\quad \phi \psi \in S^{\alpha+\beta}_{\min(m,n)}(G).$$
\item
As $Y_\ell^L$ is a linear combination of derivatives $\nabla_G^\gamma$ of length $|\gamma|=\omega_\ell$, one has (if $n\geq\omega_\ell$):
$$\phi\in S^\alpha_n(G) \quad\Longrightarrow\quad
Y_\ell^L \phi\in S^{(\alpha-\omega_\ell)_+}_{n-\omega_\ell}(G).$$
\item
As smooth functions are locally bounded, one has also:
$S^\alpha_{\alpha-1}(G) \cap C^\infty(G) \subset S^\alpha(G)$.
\end{enumerate}

\noindent
The coordinates and the coefficients of the left-invariant vector fields belong to the following classes.
\begin{prop}
One has
\begin{equation}\label{SYMB:x}
\forall \ell\in\{1,\ldots,q\},\qquad x_\ell(g) \in S^{\omega_\ell}(G)
\end{equation}
and
\begin{equation} \label{SYMB:zeta} 
\forall \ell,\ell'\in\{1,\ldots,q\},\qquad  \zeta_{\ell,\ell'}(x_1(g),\ldots,x_q(g)) \in S^{\omega_{\ell'}-\omega_\ell}(G).
 \end{equation}
\end{prop}
\paragraph{Proof.}
We already observed that $|x_\ell(g)|\lesssim \norme[G]{g}^{\omega_\ell}$
thus $x_\ell(g)\in S^{\omega_\ell}_0(G)$. Next, using~\eqref{HDERIVATIVE}, one gets:
\[
Y^L_{\ell_0}(x_\ell(g)) = \delta_{\ell_0,\ell}+\zeta_{\ell_0,\ell}(x(g))
= \begin{cases}
\zeta_{\ell_0,\ell}(x(g)) & \text{if }\omega_{\ell_0}<\omega_\ell, \\
\delta_{\ell_0,\ell} & \text{if }\omega_{\ell_0}=\omega_\ell,\\
0 & \text{if }\omega_{\ell_0}>\omega_\ell.
\end{cases}
\]
Assuming $\ell_0\in\{1,\ldots,n_1\}$, one gets $|\nabla_G x_\ell(g)| \leq C \|g\|_G^{\omega_\ell-1}$ thus $x_\ell(g)\in S^{\omega_\ell}_1(G)$.
One can now bootstrap this partial result in the expression~\eqref{ZETA}, which gets us $\zeta_{\ell,\ell'}(x(g)) \in S^{\omega_{\ell'}-\omega_\ell}_1(G)$.
The previous expression now reads $\nabla_G x_\ell(g) \in S^{\omega_\ell -1}_1$ and thus $x_\ell(g)\in S^{\omega_\ell}_2(G)$.
Iterating this process leads to $x_\ell(g)\in S^{\omega_\ell}_\infty(G)$ and $\zeta_{\ell,\ell'}(x(g)) \in S^{\omega_{\ell'}-\omega_\ell}_\infty(G)$.
\cqfd

\bigskip\par\noindent
The key result is that one can adjust the gauge norm to be a symbol of order~1 (see also \S\ref{S:vectFields}).
\begin{prop}\label{prop_symbol}
There exists a symbol $\rho(g) \in S^1_1(G)$ that is uniformly equivalent to the gauge norm. For higher-order derivatives, it satisfies
for any  multi-index $\gamma$:
\begin{equation}\label{HIGHDERRHO}
\forall g\in G, \qquad \rho(g)\leq 1 \qquad\Longrightarrow\qquad
|\nabla_G^\gamma \rho| \leq \frac{C_\gamma}{\rho^{|\gamma|-1}}\cdotp
\end{equation}
Moreover, there exists $w\in\N$ such that $\rho^w \in S^w(G)$.
\end{prop}

\paragraph{Proof.}
Let us now modify the gauge norm \eqref{GAUGE1} into the uniformly equivalent gauge:
\begin{equation}\label{GAUGE2}
\rho(g)
=\left(\sum_{\ell=1}^q |x_\ell(g)|^{w/\omega_\ell}\right)^{1/w}
\end{equation}
with $w=2\operatorname{LCM}(1,\ldots,m)$ to ensure that each $w/\omega_\ell\in 2\mathbb{N}$.
In particular, $\rho(g)^w \in C^\infty(G)$.
Next, one computes the first horizontal derivative of the norm, using~\eqref{HDERIVATIVE}:
$$\nabla_G \rho=\frac{\nabla_G(\rho^w)}{w\rho^{w-1}}
=\left(\frac{1}{\rho^{w-1}}\left[ x_\ell^{w-1}+\sum_{\ell'=1}^q\frac{\zeta_{\ell,\ell'}(x)\,x_{\ell'}^{\frac{w}{\omega_{\ell'}}-1}}{\omega_{\ell'}}\right]\right)_{\ell=1\ldots n_1}.$$
The expression in square brackets is a symbol of class $S^{w-1}(G)$ because of \eqref{SYMB:x} and \eqref{SYMB:zeta}
and $\omega_\ell=1$ for~$\ell\leq n_1$.
Thus $\nabla_G \rho$ is bounded near the origin which means that
the modified gauge $\rho$ belongs to~$S^1_1(G)$.
Next, one observes that for any $\alpha\geq1$, if $\theta\in S^\alpha(G)$ then
\[
\nabla_G\left(\frac{\theta}{\rho^\alpha}\right) = \left( \frac{\nabla_G \theta}{\rho^{\alpha-1}} -\alpha \frac{\theta}{\rho^\alpha} \nabla_G \rho \right) \frac{1}{\rho}
= \left(  \frac{\theta_1}{\rho^{\alpha-1}} + \frac{\theta_2}{\rho^{\alpha+w-1}} \right) \frac{1}{\rho}
\]
with $\theta_1\in S^{\alpha-1}(G)$ and $\theta_2\in S^{\alpha+w-1}(G)$.
One can thus claim by recurrence on the length of the multi-index $\gamma$ that
\[
\nabla_G^\gamma \rho = \left(\theta_{\gamma,0}+ \sum_k \frac{\theta_{\gamma,k}}{\rho^{\alpha_k}} \right) \frac{1}{\rho^{|\gamma|-1}}
\]
where $\theta_{\gamma,k}\in S^{\alpha_k}(G)$ is a polynomial in $x_\ell(g)$ with $\alpha_k\geq1$ and $\theta_{\gamma,0}$ is a polynomial.
Note that a polynomial in $S^0(G)$ is necessarily the sum of a constant and a polynomial in $S^1(G)$ and that, by~\eqref{HDERIVATIVE}, the horizontal
derivatives of a polynomial are also a polynomial.
This ensures~\eqref{HIGHDERRHO}. 
The final assertion about $\norme[G]{\cdot}^w$ follows immediately from~\eqref{SYMB:x}.
\cqfd

\bigskip\noindent%
From now one, one will modify the gauge norm accordingly 
and assume that $\norme[G]{\cdot} = \rho \in S^1_1(G)$.

\subsection{The radial vector field}

The \textsl{infinitesimal generator of dilations} on $\mathfrak{g}$ is the linear operator $\tilde{R}:\mathfrak{g}\to\mathfrak{g}$ defined by
\begin{equation}
\tilde{R} =\sum_{j=1}^m j\pi_j.
\end{equation}
It is diagonalizable with positive eigenvalues; its trace $\operatorname{Tr}\tilde{R}=Q$ is the homogeneous dimension. 
One checks immediately that $\delta_r  = e^{(\log r)\tilde{R}}$ thus
$\tilde{R}(x)=\left.\frac{d}{dr}\delta_r(x)\right|_{r=1}.$
The pair $(x,\tilde{R}(x))$ is a vector field on~$\mathfrak{g}$ whose expression in exponential coordinates follows from~\eqref{PROJECTIONS}:
\begin{equation}\label{RADIALEXPCOORD}
\forall x=(x_1,\ldots,x_q)\in\mathfrak{g},\qquad
\tilde{R}(x) = \sum_{\ell=1}^q \omega_\ell x_\ell \partial_\ell.
\end{equation}
Its exponential lift is called the \textsl{radial field} on $G$:
\begin{equation}
R(g)= d\exp_{\vert v}(\tilde{R}(v))=\left.\frac{d}{dr}(r g)\right|_{r=1} .
\end{equation}
\begin{prop}\label{P:RADIALFIELD}
The radial vector field is scaling invariant:
\begin{equation}\label{RSCALING}
R(rg)=\left((L_{rg})_\ast\circ\,\delta_r\circ(L_g)_\ast^{-1}\right)\left(R(g)\right).
\end{equation}
Moreover, it can be expressed in terms of left-invariant derivatives:
\begin{equation}\label{RCOORD}
R(g)=\sum_{\ell=1}^q \sigma_\ell(x_1(g),\ldots,x_q(g))Y_\ell^L(g)
\end{equation}
with
$\displaystyle
\sigma_\ell(x(g))=
\omega_\ell x_\ell(g)
+
\sum_{n=1}^{m-1}\frac{(-1)^n}{(n+1)!}
\sum_{\ell_1,\ldots,\ell_{n+1}}
x_{\ell_1}(g) \cdots x_{\ell_{n+1}}(g)
\cdot
\omega_{\ell_{n+1}} \kappa(\ell_1,\ldots,\ell_{n+1};\ell)
\in S^{\omega_\ell}(G)$.
\end{prop}

\paragraph{Remarks}
\begin{itemize}
\item
Note that the variable $x_{\ell_i}$ that appears in the second term defining $\sigma_\ell$ must satisfy
\[
\omega_\ell = \omega_{\ell_1}+\ldots+\omega_{\ell_{n+1}}
\]
because if it is not the case, then $\kappa(\ell_1,\ldots,\ell_{n+1};\ell)=0$. In particular, as there are at least $n+1\geq 2$ factors, one has $\omega_{\ell_i}<\omega_\ell$
for each $i$.
\item Both expressions for $\tilde{R} = \sum \omega_\ell x_\ell \partial_\ell = \sum \sigma_\ell \tilde{Y}_\ell^L$
combined with~\eqref{HDERIVATIVE} provide a remarkable identity embedded in the commutator structure. For any $\ell'\in\{1,\ldots,q\}$:
\begin{gather*}
\sum_{n=1}^{m-1}\frac{(-1)^n}{(n+1)!}
\sum_{\ell_1,\ldots,\ell_{n+1}}
x_{\ell_1} \cdots x_{\ell_{n+1}} \omega_{\ell_{n+1}}
\left(
\kappa(\ell_1,\ldots,\ell_{n+1};\ell') + \sum_{\ell} \kappa(\ell_1,\ldots,\ell_{n+1};\ell) \zeta_{\ell,\ell'}(x)
\right) \\
= -\sum_{\ell} \omega_\ell x_\ell \zeta_{\ell,\ell'}(x).
\end{gather*}
Note that when one substitutes $x=x(g)$, both sides are indeed symbols of class $S^{\omega_{\ell'}}(G)$.
\end{itemize}
\par

\paragraph{Proof.}
Formula~\eqref{RSCALING} follows~\textit{e.g.}
from the identities $\tilde{R}\circ \delta_r = \delta_r \circ \tilde{R}$ and $\operatorname{ad}\circ\,\delta_r=\delta_r\circ\operatorname{ad}\circ\,\delta_r^{-1}$:
\begin{align*}
R(rg) &= d\exp_{\vert \delta_r(v)}\circ\,\delta_r\circ \tilde{R} (v)\\
&= (L_{rg})_\ast \circ \left(\frac{1-e^{-\operatorname{ad}(\delta_r(v))}}{\operatorname{ad}(\delta_r(v))}\right)\circ\,\delta_r\circ \tilde{R} (v)\\
&=(L_{rg})_\ast \circ\,\delta_r\circ \left(\frac{1-e^{-\operatorname{ad}(v)}}{\operatorname{ad}(v)}\right)\circ \tilde{R} (v)\\
&=\left((L_{rg})_\ast\circ\,\delta_r\circ(L_g)_\ast^{-1}\right)\left(R(g)\right).
\end{align*}
The definition of $R(g)$ with $g=\exp(v)$ also reads
$$R(g)= \sum_{\ell=1}^q \omega_\ell x_\ell(g) (d\exp_{\vert v}Y_\ell).$$
Combining the expression for the differential of $\exp$ given by~\eqref{DIFFEXP},
the identity $[u,v]^L=[u^L,v^L]$ and the fact that $v=\sum x_\ell(g) Y_\ell$ give:
\begin{align*}
R(g) &=\sum_{\ell=1}^q  \omega_\ell x_\ell(g) \left(\sum_{n=0}^{m-1}
\frac{ (-1)^n}{(n+1)!} \underset{n\text{ times }v}{[v,\ldots,[v},Y_\ell]]\right)^L\!\!\!(g)\\
&=\sum_{\ell=1}^q  \omega_\ell x_\ell(g) \left(
Y_\ell^L(g) +
\sum_{n=1}^{m-1}\sum_{\ell_1,\ldots,\ell_n}
\frac{ (-1)^n}{(n+1)!} x_{\ell_1}(g)\ldots x_{\ell_n}(g) [Y_{\ell_{1}},\ldots,[Y_{\ell_n},Y_\ell]]^L(g)
\right).
\end{align*}
This formula can be further simplified into \eqref{RCOORD} by using \eqref{COMMSTRUCT}.
The symbol property comes from \eqref{SYMB:x} and the restriction on non-vanishing indices imposed by~\eqref{COMMSTRUCT}.
 \cqfd\par

\begin{exs}
The previous computation can be simplified further by observing the anti-symmetrical role of $\ell_n$ and $\ell_{n+1}$ in $\omega_{\ell_{n+1}}[Y_{\ell_n},Y_{\ell_{n+1}}]$
if $\omega_{\ell_n}=\omega_{\ell_{n+1}}$. For $m\leq 4$, one thus gets the following expressions for the radial field.
\begin{enumerate}
\item
For a group of step $m=2$, the radial field is given by:
$$R(g)
=\sum_{\ell=1}^q \omega_\ell x_\ell(g) Y_\ell^L(g).$$
On the Heisenberg group $\HH$ with exponential coordinates introduced in \S\ref{S:expCoord}, this formula boils down, as expected, to the following one:
\[
R(g)=\left(pY_1+qY_2+2(r-\frac{1}{2}pq)Y_3\right)^L\!\!\!(g) =
p \partial_p + q \partial_q + 2r \partial_r.
\]
\item
For a group of step $m=3$, the radial field is ``corrected'' along $V_3$:
\[
R(g)=\sum_{\ell=1}^q \omega_\ell x_\ell(g) Y_\ell^L -\frac{1}{2}\sum_{\substack{1\leq \ell_1\leq n_1\\ n_1<\ell_2\leq n_2}} x_{\ell_1}(g)x_{\ell_2}(g)[Y_{\ell_1},Y_{\ell_2}]^L.
\]
\item For step $m=4$, its expression involves a further ``correction'' along $V_4$ that is split among two types of commutators:
\begin{align*}
R(g)= &\sum_{\ell=1}^q \omega_\ell x_\ell(g) Y_\ell^L -\frac{1}{2}\sum_{\substack{1\leq \ell_1\leq n_1\\ n_1<\ell_2\leq \mathbf{n_3}}} x_{\ell_1}(g)x_{\ell_2}(g)[Y_{\ell_1},Y_{\ell_2}]^L\\&+\frac{1}{6}\sum_{\substack{1\leq \ell_1,\ell_2\leq n_1\\ n_1<\ell_3\leq n_2}} x_{\ell_1}(g)x_{\ell_2}(g)x_{\ell_3}(g)[Y_{\ell_1},[Y_{\ell_2},Y_{\ell_3}]]^L.
\end{align*}
\end{enumerate}
\end{exs}

\begin{prop}
The gauge norm \eqref{GAUGE2} and the radial field are related by the following formula:
\begin{equation}\label{RIPP}
\forall s>0,\qquad 
\frac{1}{\norme[G]{\cdot}^{2s}} = -\frac{1}{2s} R\left(\frac{1}{\norme[G]{\cdot}^{2s}}\right).
\end{equation}
\end{prop}
\paragraph{Proof.}
Applying the chain rule, one gets:
\[
\frac{\lambda(g)}{\norme[G]{g}^{2s}} = -\frac{1}{2s} R\left(\frac{1}{\norme[G]{g}^{2s}}\right)
\]
with $\lambda(g)=\frac{R(\norme[G]{g})}{\norme[G]{g}}$ and where the field $R$ is obviously computed at the same point $g\in G$
as the function that is being derivated.
Let us also observe that:
\[
\lambda(g)=\frac{R\left(\norme[G]{g}^w\right)}{w\norme[G]{g}^w}
\]
for any $w\in\mathbb{N^\ast}$ and in particular for $w=2\operatorname{LCM}(1,\ldots,m)$ for which we know that $\norme[G]{\cdot}^w \in S^w(G)$
by proposition~\ref{prop_symbol}. 
Using the formula~\eqref{GAUGE2} for the modified gauge norm and~\eqref{RADIALEXPCOORD} for the expression of the radial field
in exponential coordinates, one then gets (note that $w/\omega_\ell\in 2\N^\ast$):
\[
\forall x\in\mathfrak{g}\simeq\R^q,\qquad
\tilde{R}(\norme[\mathfrak{g}]{x}^w)
= \sum_{\ell=1}^q \omega_\ell x_\ell \cdot  \frac{w}{\omega_\ell} x_\ell^{\frac{w}{\omega_\ell}-1}
= w\norme[\mathfrak{g}]{x}^w
\]
and thus $\lambda(g)= 1$ for any $g\in G$.
 \cqfd\par

\subsection{Adjoints}

\begin{prop}\label{ADJOINTR}
For the $L^2(G)$ scalar product, the adjoint vector field to $R$ is
\[
R^\ast(g) = -Q -R(g).
\]
\end{prop}
\paragraph{Proof.}
The proof is simplest in exponential coordinates, using~\eqref{RADIALEXPCOORD} and  \eqref{HOMDIM}:
\[
\forall x\in\mathfrak{g}\simeq\R^q,\qquad
\tilde{R}(x)+\tilde{R}^\ast(x) = \div\tilde{R} = \sum_{\ell=1}^q \omega_\ell = Q.
\]
One can also prove this formula directly, using~\eqref{RCOORD} and \eqref{HDERIVATIVE}:
\[
\forall g\in\mathfrak{G},\qquad
R(g)+R^\ast(g) = \sum_{\ell} \partial_\ell(\sigma_\ell) + \sum_{\ell,\ell'} \zeta_{\ell,\ell'} \cdot (\partial_{\ell'} \sigma_\ell)
+ \sigma_\ell \cdot (\partial_{\ell'} \zeta_{\ell,\ell'}).
\]
In this sum, according to a remark that follows~\eqref{HDERIVATIVE},
the index $\ell'$ is restricted to $\ell'>n_{\omega_\ell}$ and, in particular, the definition \eqref{def:njomegaj} then implies $\omega_{\ell'}>\omega_\ell$.
Now thanks to the remark that follows proposition~\ref{P:RADIALFIELD}, one can claim that the variable $x_\ell$ does not appear in the second part of $\sigma_\ell$, thus
its derivative reads
\[
\partial_\ell (\sigma_\ell) = \omega_\ell.
\]
For a similar reason, $\partial_{\ell'} \sigma_\ell =0$ for $\omega_{\ell'}>\omega_\ell$.
One observes also that in~\eqref{ZETA}, each $\ell_i$ involved in the expression of $\zeta_{\ell,\ell'}$ must satisfy $\omega_{\ell_i}<\omega_{\ell'}$.
In particular, $\partial_{\ell'} \zeta_{\ell,\ell'}=0$.
One concludes using~\eqref{HOMDIM}.
 \cqfd\par

\bigskip
The next property checks that left-invariant vector fields on a Carnot group are divergence-free.
\begin{prop}\label{PROPADJFIELD}
For the $L^2(G)$ scalar product, the adjoint vector field to $Y_\ell^L$ is $-Y_\ell^L$.
In particular, for any smooth function $\psi$ on $G$ and any $\ell_1,\ldots,\ell_n\in\{1,\ldots,q\}$:
\begin{equation}
\int_{G} [Y_{\ell_1}^L,\ldots,[Y_{\ell_{n-1}}^L,Y_{\ell_n}^L]] \psi(g) \cdot \psi(g) dg= 0.
\end{equation}
\end{prop}
\paragraph{Proof.}
The second ``computational'' proof of the previous proposition (the one based on \eqref{HDERIVATIVE}) also ensures that
\[
\partial_{\ell'} \zeta_{\ell,\ell'}=0
\]
when $\omega_{\ell'}>\omega_\ell$ and therefore $(Y_\ell^L)^\ast = - Y_\ell^L$.
As the commutator of two antisymmetric operators is also an antisymmetric one, the second statement follows immediately.
\cqfd\par

\subsection{A density result}

The following density result can be proved by a scaling argument.

\begin{prop}
The space $\mathcal{D}(G\backslash\{e\})$ of $C^\infty$ functions, compactly supported outside the origin,
is dense in $H^s(G)$ for any $0\leq s<Q/2$.
\end{prop}

\paragraph{Proof.}
One can use a Hilbert space approach based on scaling and Schwartz's theorem for distributions.
Let us assume additionally that $s\in\N$ and consider a function $u\in H^s(G)$ that is
orthogonal to any $\phi\in \mathcal{D}(G\backslash\{e\})$, \textit{i.e.}:
$$(u,\phi)_{s}=\sum_{|\gamma|\leq s} \int_G \nabla_G^\gamma u(g) \cdot \nabla_G^\gamma\phi(g) dg=0.$$
Integrating by parts (using proposition~\ref{PROPADJFIELD} and the notation $\gamma^\ast$ for the multi-index $\gamma$ in reverse order) reads:
$$\sum_{|\gamma|\leq s} (-1)^{|\gamma|} \int_G \nabla_G^{\gamma^\ast} \nabla_G^\gamma u(g) \cdot \phi(g) dg=0.$$
For fractional values of $s$, one would replace $\nabla_G^{\gamma^\ast}\nabla_G^\gamma$ by a fractional power of the sub-Laplacian~\eqref{HYPOLAPLACE}
and what follows would go unchanged.
Schwartz's theorem implies that the distributional support of
$$v=\sum_{|\gamma|\leq s} (-1)^{|\gamma|} \nabla_G^{\gamma^\ast}\nabla_G^\gamma u$$
is reduced to the single point $\{e\}$ and thus $v=\sum (-1)^{|\alpha|} c_\alpha \partial^\alpha\delta$ where $\delta$ is the Dirac function at the origin.
As $v$ is at most a $2s^{\text{th}}$ horizontal derivative of $u$, one has $v\in H^{-s}(G)$  and in particular for any test function $\psi\in\mathcal{D}(G)$:
$$\left|\int_G v(g)\psi(g) dg\right|^2 \leq C \sum_{|\gamma|\leq s} \int_G |\nabla_G^\gamma \psi(g)|^2 dg.$$
The constant $C$ does not depend on the support of $\psi$ because $\operatorname{supp} v \subset \{e\}$.
In particular, one can apply this inequality to the dilations $\psi(r g)$ for any $r>1$:
$$\left|\int_G v(g)\psi(r g) dg\right|^2  \leq C \sum_{|\gamma|\leq s} r^{2|\gamma|} \int_G |\nabla_G^\gamma \psi(r g)|^2 dg.$$
thus
$$\left|\int_G v(r^{-1} g)\psi(g) dg\right|^2  \leq C\sum_{|\gamma|\leq s} r^{2|\gamma|+Q} \int_G |\nabla_G^\gamma \psi( g)|^2 dg.$$
Finally, one can compute the left-hand side using the homogeneity of the Dirac mass:
$$\int_G v(r^{-1} g)\psi(g) dg = \sum_\alpha c_\alpha r^{Q+\sum \alpha_j\omega_j}  \partial^\alpha\psi(e).$$
Combining both formulas, one gets for any $r>1$:
$$\left|\sum_\alpha c_\alpha r^{Q+\sum \alpha_j\omega_j}  \partial^\alpha\psi(e)\right| \leq C r^{Q/2+s}\norme[H^s(G)]{\psi}$$
and in particular with a suitable choice of $\psi$ and $r\to\infty$:
$$c_\alpha\neq0 \qquad\Longrightarrow\qquad s \geq \frac{Q}{2}+\sum \omega_j\alpha_j.$$
But as $s<Q/2$, each coefficient $c_\alpha$ vanishes, \ie $v=0$ in $H^{-s}(G)$ and thus using $u\in H^s(G)$ as a test function, one infers $u=0$.
\cqfd

\begin{rmk}
When $Q$ is even and $s=Q/2\in\mathbb{N}$, the previous density result still holds. The only change in the proof is to observe that $\delta\not\in H^{-Q/2}(G)$ by exhibiting an example of an unbounded function
in $H^{Q/2}(G)$; the classical example $\log(-\log \norme{g})\psi(g)$ with a sooth cut-off $\psi$ still works. However, when $Q$ is odd, one still has $\delta\not\in H^{-Q/2}(G)$ but the density result \textbf{fails} as it already does in $H^{n+\frac{1}{2}}(\R^{2n+1})$. For more details on this point, see~\cite{PhD}.
\end{rmk}

\subsection{Hardy inequality}\label{PROOF}

In this final section, let us combine the previous results into a proof of Theorem~\ref{main_thm}.

\bigskip
Given $f\in H^s(G)$ with $s<Q/2$ and the density result of the previous section,
one can assume without restriction that $f$ is compactly supported and that $0\notin\operatorname{Supp}u$.
Next, one will take a smooth cutoff function $\chi:\R\to[0,1]$ such that $\chi(t)=1$ if $|t|<1/2$. For any
$\rho_0>0$, one has:
\begin{equation}\label{SUPPRESTR}
\int_G \frac{|f(g)|^2}{\norme[G]{g}^{2s}} \leq \int_G \frac{|\varphi(g)|^2}{\norme[G]{g}^{2s}}  + 
\left(\frac{2}{\rho_0}\right)^{2s} \norme[L^2(G)]{f}^2
\quad\text{with}\quad \varphi(g)=\chi\left(\frac{\norme[G]{g}}{\rho_0}\right)f(g).
\end{equation}
Moreover, one has $\norme[H^s(G)]{\varphi}\leq C_s\rho_0^{-s}\norme[H^s(G)]{f}$. Without restriction, one
can therefore assume that $f$ (now denoted by $\varphi$) is compactly supported in a fixed but arbitrary small annular neighborhood around the
origin.

\bigskip
The key of the computation is the following integration by part argument.
Using~\eqref{RIPP}, one has
\[
\int_G \frac{|\varphi(g)|^2}{\norme[G]{g}^{2s}}
= -\frac{1}{2s}\int_G R\left(\frac{1}{\norme{g}^{2s}}\right) \cdot |\varphi(g)|^2.
\]
Using proposition \ref{ADJOINTR} and the fact that $\operatorname{supp}\phi$ is an annulus around the origin so that no boundary terms appear:
\[
\left(\frac{Q}{2}-s\right)\int_G \frac{|\varphi(g)|^2}{\norme[G]{g}^{2s}}
= -\int_G \frac{\varphi(g) R(\varphi(g))}{\norme{g}^{2s}}\cdotp
\]
According to  \eqref{RCOORD}, the radial field can be expressed with left-invariant vector fields:
\begin{equation}\label{KEYIBP}
\left(\frac{Q}{2}-s\right)\int_G \frac{|\varphi(g)|^2}{\norme[G]{g}^{2s}}
= - \sum_{\ell=1}^q
\int_G \frac{\sigma_\ell(x(g)) \varphi(g) Y_\ell^L(\varphi(g))}{\norme{g}^{2s}}\cdotp
\end{equation}
What we do next depends on the order of each derivative $Y_\ell^L\simeq \nabla_G^{\omega_\ell}$.

\paragraph{Case $m=1$.}
In the Euclidian case, one uses Cauchy-Schwarz and Young's identity $|ab| \leq \varepsilon a^2 + \varepsilon^{-1} b^2$ with $\varepsilon>0$ small enough
so that $s+\varepsilon<Q/2$, which leads to:
\begin{equation}\label{CASEm1}
\left(\frac{Q}{2}-s-\varepsilon \right)\int_G \frac{|\varphi(g)|^2}{\norme[G]{g}^{2s}}
\leq \varepsilon^{-1} C \int_G \frac{|\nabla_G \varphi|^2}{\norme[G]{g}^{2(s-1)}}.
\end{equation}
This proves Hardy's inequality for $s=1$. Interpolation with $L^2$ then ensures that the Hardy inequality holds for any $s\in[0,1]$.
Finally, the previous estimate provides a bootstrap argument from $s-1$ to $s$ for any $s<Q/2$.

\paragraph{Case $m= 2$.}
One uses the Euclidian technique to deal with the horizontal derivatives. For the stratum~$V_2$, one uses the commutator structure
to backtrack one ``half'' integration by part. More precisely, the right-hand side of \eqref{KEYIBP} becomes for $1\leq \ell \leq n_1$:
\[
\left|
\int_G \frac{\sigma_\ell(x(g)) \varphi(g) Y_\ell^L(\varphi(g))}{\norme{g}^{2s}}
\right| \leq 
\varepsilon \int_G \frac{|\varphi(g)|^2}{\norme[G]{g}^{2s}}
+ C_\varepsilon \int_G \frac{|\nabla_G \varphi(g)|^2}{\norme[G]{g}^{2(s-1)}}
\]
and  using~\eqref{COMMREVSTRUCT} for $n_1<\ell\leq n_2=q$ and proposition~\ref{PROPADJFIELD} (notice the cancellation of the highest order term
thanks to the commutator structure):
\begin{align*}
-\int_G \frac{\sigma_\ell(x(g))}{\norme{g}^{2s}}\:
\varphi(g) \cdot [Y_{\alpha_1(\ell)}^L,Y_{\alpha_2(\ell)}^L](\varphi(g))
=
& \int_G  Y_{\alpha_1(\ell)}^L\left( \frac{\sigma_\ell(x(g))}{\norme{g}^{2s}}\right) \:
 \varphi(g) \cdot Y_{\alpha_2(\ell)}^L (\varphi(g)) \\
& -\int_G  Y_{\alpha_2(\ell)}^L\left( \frac{\sigma_\ell(x(g))}{\norme{g}^{2s}}\right) \:
 \varphi(g) \cdot Y_{\alpha_1(\ell)}^L (\varphi(g)).
\end{align*}
Using the symbol properties of $\norme[G]{\cdot}$ and $\sigma_\ell(x(g))$, both terms are bounded in the following way:
\[
\int_G \frac{ |\varphi(g)| \cdot |\nabla_G\varphi(g) |}{\norme[G]{g}^{2s-1}}
\leq 
\varepsilon \int_G \frac{|\varphi(g)|^2}{\norme[G]{g}^{2s}}
+ C_\varepsilon \int_G \frac{|\nabla_G \varphi(g)|^2}{\norme[G]{g}^{2(s-1)}}\cdotp
\]
One thus gets~\eqref{CASEm1} again and one can conclude the proof just as in the case $m=1$.

\paragraph{Case $m\geq 3$.}
The additional terms on the right-hand side of~\eqref{KEYIBP} correspond to $n_2< \ell \leq q$. Thanks to~\eqref{COMMREVSTRUCT},
one can express each of them with commutators from the first stratum:
\[
I_\ell(\varphi) = -\int_G \frac{\sigma_\ell(x(g))}{\norme{g}^{2s}}\:
\varphi(g) \cdot [Y_{\alpha_1(\ell)}^L,\ldots,[Y_{\alpha_{\omega_\ell-1}(\ell)}^L,Y_{\alpha_{\omega_\ell}(\ell)}^L]] (\varphi(g)).
\]
As in the case $m=2$, the key is to use the commutator structure to put all the derivatives but one on the symbol.
More precisely, using proposition~\ref{PROPADJFIELD}, one first gets:
\[
I_\ell(\varphi) = \frac{1}{2}\int_G 
 [Y_{\alpha_1(\ell)}^L,\ldots,[Y_{\alpha_{\omega_\ell-1}(\ell)}^L,Y_{\alpha_{\omega_\ell}(\ell)}^L]] 
 \left(
 \frac{\sigma_\ell(x(g))}{\norme{g}^{2s}}
 \right)\:
|\varphi(g)|^2.
\]
Next, one puts the outermost derivative back out onto $\varphi(g)^2$:
\[
I_\ell(\varphi) = -\sum_{i=1}^{\omega_\ell} \int_G 
W_{\ell,i}
 \left(
 \frac{\sigma_\ell(x(g))}{\norme{g}^{2s}}
 \right)\:
\varphi(g) \cdot Y_{\alpha_i(\ell)}^L \varphi(g)
\]
where each $W_{\ell,i}$ is a derivative of order $\omega_{\ell}-1$. The symbol property $\sigma_\ell(x(g))\in S^{\omega_\ell}(G)$
given by proposition~\ref{P:RADIALFIELD} ensures that
\[
\left| W_{\ell,i} \left(  \frac{\sigma_\ell(x(g))}{\norme{g}^{2s}} \right) \right|
\leq \frac{C_{\ell,i,s}}{\norme[G]{g}^{2s-1}}\cdotp
\]
Again, one gets
\[
|I_\ell(\varphi)| \leq 
\varepsilon \int_G \frac{|\varphi(g)|^2}{\norme[G]{g}^{2s}}
+ C_\varepsilon \int_G \frac{|\nabla_G \varphi(g)|^2}{\norme[G]{g}^{2(s-1)}}
\]
and~\eqref{CASEm1} holds once more.
As in the case $m=1$, one thus gets the Hardy inequality for $s=1$. Then, by interpolation with $L^2(G)$, one gets it for $s\in[0,1]$.
Finally, using~\eqref{CASEm1} iteratively, one can collect it for any $s<Q/2$.

\begin{rmk}
When $Q$ is odd and $s\in\N$, one can always take $\varepsilon$ small enough so that $\frac{Q}{2}-s-\varepsilon\neq 0$ in~\eqref{CASEm1}.
The previous iteration argument thus proves the Hardy inequality \eqref{MAINHARDY} for any $s\in\N$, but the result is then
only valid for functions that belong to the $H^s(G)$-closure
of smooth compactly supported functions whose support avoids the origin.
\end{rmk}

\subsection{Homogeneous Hardy inequality}

One can slightly improve~\eqref{MAINHARDY} by using a simple scaling argument.
For simplicity, we will only spell out the procedure for $s\in\N$ and $0\leq s<Q/2$ though it would also work for fractional
values of $s$ if $\nabla_G^s$ was replaced by the corresponding power $\mathcal{L}_G^{s/2}$ of the subellliptic Laplace operator~\eqref{HYPOLAPLACE}.
\begin{thm}\label{hom_thm}
For $0\leq s<Q/2$, the following homogeneous inequality holds:
\begin{equation}\label{HOMOGENEOUSHARDY}
\forall f\in H^s(G)\qquad
\int_G \frac{|f(g)|^2}{\norme{g}^{2s}} dg \leq 2C_s \norme[L^2(G)]{\nabla_G^s f}^2.
\end{equation}
\end{thm}

\paragraph{Proof.}
Let us indeed apply~\eqref{MAINHARDY} to the function $f(r^{-1}g)$. After the change of variable $g=r\bar{g}$, one gets:
\[
r^{Q-2s} \int_G \frac{|f(\bar{g})|^2}{\norme{\bar{g}}^{2s}} d\bar{g}
\leq C_s
\sum_{|\alpha|\leq s} r^{Q-2|\alpha|} \int_G |\nabla_G^\alpha f(\bar{g})|^2 d\bar{g}
\]
which, for $r<1$, can be further simplified into:
\[
\int_G \frac{|f(\bar{g})|^2}{\norme{\bar{g}}^{2s}} d\bar{g}
\leq C_s
\sum_{|\alpha|=s}\int_G |\nabla_G^\alpha f(\bar{g})|^2 d\bar{g} 
+
C_s r^2
\sum_{|\alpha|\leq s-1} \int_G |\nabla_G^\alpha f(\bar{g})|^2 d\bar{g}.
\]
Choosing
\[
r^2 = \min\left\{ 1 ; \frac{\norme[L^2(G)]{\nabla_G^s f}^2}{\norme[H^{s-1}(G)]{f}^2}
\right\}
\]
instantly leads to~\eqref{HOMOGENEOUSHARDY}.
\cqfd\par

\begin{rmk} It would have been tempting to try using~\eqref{SUPPRESTR}-\eqref{KEYIBP} without digging further in the commutator structure to get
\[
\int_G \frac{|f(g)|^2}{\norme[G]{g}^{2s}} dg
\leq  \int_G |f(g)|^2 dg+
\left(\frac{Q}{2}-s -\varepsilon\right)^{-1} \sum_{\ell=1}^q
\int_G \frac{| Y_\ell^L((\chi f)(g))|^2 }{\norme{g}^{2(s-\omega_\ell)}} dg.
\]
For $s=1$, it gives a Hardy inequality with $\norme[H^m(G)]{f}^2$ on the right-hand side.
However, a scaling argument is then not sufficient to deduce the correct one, either~\eqref{MAINHARDY} or \eqref{HOMOGENEOUSHARDY}.
Indeed, one would simultaneously need to let $r\to\infty$ and $r\to0$ to get rid of the superfluous derivatives without letting the lower-order $L^2(G)$
term get in the way, which is overall impossible. 
\end{rmk}

\section{A remark about the case of general hypoelliptic vector fields}\label{S:vectFields}

For general families of vector fields that satisfy a H\"ormander condition of step $m$, the technique of proving
the Hardy inequality by integration by part works, but possibly with some restrictions.

\subsection{A counter-example to the symbol property of the gauge}

The main objection is the following one. When the group structure is discarded, the fact that one can chose
a gauge pseudo-norm in a symbol class of order 1 can fail.

\bigskip
For example, the family
\[
Z_1=\partial_1+x_1\partial_3, \quad
Z_2=\partial_2+x_4\partial_3+x_5 \partial_4
\quad\text{and}\quad
Z_3=\partial_5
\]
is uniformly of rank $3$ in $\R^5$ and satisfies a uniform H\"ormander bracket condition of step~3:
\[
\partial_4=[Z_3,Z_2], \qquad \partial_3=[[Z_3,Z_2],Z_2].
\]
However,
the ``natural'' gauge:
$$\rho=(|x_1|^{12}+|x_2|^{12}+|x_3|^4+|x_4|^6+|x_5|^{12})^{1/12}$$
is \textbf{not} a symbol of order 1 because $|Z_1\rho|\geq c\rho^{-1}$ along $x_1^3-x_3=x_2=x_4=x_5=0$.
Luckily, for this particular family, the change of variable $y_3=x_3-\frac{1}{2}x_1^2$ and $y_i=x_i$ ($i\neq3$)
transforms the family into $Z_1'=\partial_{y_1}$, $Z'_2=\partial_{y_2}+y_4\partial_{y_3}+y_5\partial_{y_4}$
and $Z'_3=\partial_{y_5}$ and for this new family, the associated gauge is a symbol of order 1.

\bigskip
In \cite[chap.~7]{PhD}, it was shown that up to a H\"ormander condition of step 3, one can
always modify the gauge by a local diffeomorphism to restore the symbol property.
However, the same question for a family of vector fields that satisfy a H\"ormander condition of step 4 or higher is still open.
For the convenience of the reader, we will recall here briefly the key points of the discussion (and clarify the redaction),
as this result was written in French and never published.

\subsection{Regular hypoelliptic vector fields of step $m$}

Let us consider a family $\mathfrak{X}=(X_\ell)_{1\leq \ell \leq n_1}$ of vector fields on some smooth open set $\Omega\subset\R^q$ and
\begin{equation}
\forall x\in\Omega,\qquad
W_k(x) = \operatorname{Span}\left(
X_i(x),\ldots, [X_{j_1},\ldots,[X_{j_{k-1}},X_{j_k}]](x)
\right).
\end{equation}
One assumes that $\underline{x}\in\Omega$ is a \textsl{regular H\"ormander point}, \ie that $n_k=\dim W_k$ is constant near $\underline{x}$
and that~$n_m=q$ for some finite integer $m\geq 2$.

\begin{rmk}
At the origin of a Carnot group \eqref{DEFSTRAT}, one would have $W_k(e)=\oplus_{j=1}^k V_j$.
\end{rmk}

Next, one introduces a local basis of vector fields $(Y_\ell(x))_{1\leq \ell\leq q}$, adapted to the stratification, \ie $Y_\ell(x) \in W_{\omega_\ell}(x)$
where for each $\ell$, the weight $\omega_\ell\in\{1,\ldots,m\}$ is defined by~\eqref{def:njomegaj}. 
For simplicity, one will now restrict $\Omega$ to be a bounded and small enough neighborhood of $\underline{x}$ on which all those properties hold.
The analog of horizontal derivatives is the family:
\begin{equation}
\nabla_{\mathfrak{X}} = (Y_1,\ldots, Y_{n_1}).
\end{equation}

A local coordinate system $(x_\ell)_{1\leq\ell\leq q}$ is said to be \textsl{adapted} to the commutator structure of the vector fields $\mathfrak{X}$ near $\underline{x}$ if
the dual basis $(\partial_\ell)_{1\leq\ell\leq q}$ satisfies $Y_\ell(\underline{x})=\partial_\ell$.

\begin{rmk}
Let us point out that \textsl{adapted} coordinates are not necessarily \textsl{privileged} in the sense of A.~Bellaiche~\cite{Bellaiche} and M.~Gromov~\cite{Gromov}:
the point of coordinates $(x_\ell)_{1\leq\ell\leq q}$ does not necessarily match with the image of $\underline{x}$ under the composite
action of the flows $e^{x_\ell Y_\ell}$ (for some predetermined order of composition).
\end{rmk}

In an adapted coordinate system, the \textsl{gauge} is defined by:
\begin{equation}
\rho(x) = \left( \sum_{\ell=1}^q |x_\ell|^{w/\omega_\ell} \right)^{1/w}
\end{equation}
where $w=2\operatorname{LCM}(1,\ldots,m)$ and the basis of vector fields and their commutators satisfy
\begin{equation}
\forall\ell\in\{1,\ldots,q\},\qquad
Y_\ell(x) = \partial_\ell + \sum_{\ell'=1}^q \zeta_{\ell,\ell'}(x) \partial_{\ell'}.
\end{equation}
One obviously has $|x_\ell|\leq \rho(x)^{\omega_\ell}$ and, using a Taylor expansion, $\zeta_{\ell,\ell'}(x_0)=0$ implies $|\zeta_{\ell,\ell'}(x)| \leq C \rho(x)$.
However, for derivatives, one can only claim that $\nabla_\mathfrak{X}^\gamma x_\ell$ and $\nabla_\mathfrak{X}^\gamma \zeta_{\ell,\ell'}$ are bounded
when $|\gamma|\geq 1$.

\subsection{A positive result for hypoelliptic fields of step 2}

\begin{thm}\label{STEP2FIELDS}
Let us consider a family of vector fields and $\underline{x}\in\Omega$ a regular H\"ormander point of step~${m=2}$.
Then for any adapted coordinate system, the gauge $\rho$ satisfies
\begin{equation}\label{FIELDSHIGHDERGAUGE}
|\nabla_\mathfrak{X}^\gamma \rho| \leq C_\gamma \rho^{1-|\gamma|}
\end{equation}
in the neighborhood of $\underline{x}$, for any multi-index $\gamma$.
\end{thm}

\paragraph{Proof.}
For $\gamma=0$, the estimate \eqref{FIELDSHIGHDERGAUGE} comes from the fact that $\rho$ is smooth and vanishes at the origin
and thus admits a Taylor expansion at the origin that is locally bounded by $\sum |x_\ell|$ and thus by $\rho$.
For~$|\gamma|=1$, the computation is actually explicit:
\[
\nabla_\mathfrak{X} \rho = \frac{1}{\rho^3} \left(
x_\ell^3 + \sum_{\ell'\leq n_1} \zeta_{\ell,\ell'}x_{\ell'}^3 + \frac{1}{2} \sum_{\ell'>n_1} \zeta_{\ell,\ell'} x_{\ell'}
\right)_{1\leq \ell\leq n_1}.
\]
In the parenthesis, the first term is locally bounded by $\rho^3$, the second by $\rho^4$ and the last one again by~$\rho^3$,
thus~$\nabla_{\mathfrak{X}}\rho\in L^\infty(\Omega)$ provided $\Omega$ is small enough.
To deal with the higher-order derivatives, let us introduce the class $\mathcal{P}_n$ of homogeneous
polynomials of $x_\ell$ and $\zeta_{\ell,\ell'}$ with smooth coefficients, \ie
\[
\sum_{\alpha,\beta} x_1^{\alpha_1} \ldots x_q^{\alpha_q} \zeta_{1,1}^{\beta_{1,1}} \zeta_{1,2}^{\beta_{1,2}} \ldots \zeta_{q,q}^{\beta_{q,q}} f_{\alpha,\beta}(x)
\]
where $f_{\alpha,\beta}\in C^\infty(\Omega)$ and $\sum \alpha_i \omega_i + \sum \beta_{j,j'} = n$.
For $n\leq 0$, one sets $\mathcal{P}_n = C^\infty(\Omega)$.
With Leibnitz formula, one checks immediately that:
\[
\partial_\ell (\mathcal{P}_n) \subset 
\begin{cases}
\mathcal{P}_{n-1} + \mathcal{P}_n & \text{if } \ell\leq n_1\\
\mathcal{P}_{n-2} + \mathcal{P}_{n-1} + \mathcal{P}_n & \text{if } \ell > n_1
\end{cases}
\]
thus $\nabla_\mathfrak{X} (\mathcal{P}_n) \subset \mathcal{P}_{n-1} + \mathcal{P}_n + \mathcal{P}_{n+1}$.
Moreover, for $m\geq n$, any expression in $\mathcal{P}_{m}$ is locally bounded by~$C\rho^n$ for some constant~$C$.
We have shown above that $\nabla_\mathfrak{X} \rho \in \rho^{-3}\cdot (\mathcal{P}_3+\mathcal{P}_4)$.
One then gets recursively on~$k=|\gamma|$ that~$\nabla_{\mathfrak{X}}^\gamma \rho$  is a linear combination of expressions
\[
\frac{\mathcal{P}_m}{\rho^{n+k-1}}
\]
with $m\geq n$ and is thus locally bounded by $C \rho^{1-k}$.
\cqfd\par

\begin{rmk}
One has $\mathcal{P}_n\subset \mathcal{P}_{n-2}$. However, for $\ell>n_1$, one has $x_\ell^2 \in \mathcal{P}_4\cap\mathcal{P}_2$ but $x_\ell^2 \not\in \mathcal{P}_3$.
\end{rmk}

\subsection{Two positive results for hypoelliptic fields of step $m\geq3$}

Let us now revert to the case of a general value for $m$.
As pointed out at the beginning of this section, one can find a counter-example of a family of vector fields, a regular H\"ormander point of step~${m=3}$
and an adapted coordinate system for which \eqref{FIELDSHIGHDERGAUGE} fails. If we tried to run the previous proof, the failure point would be that
\[
\partial_{\ell'}(\mathcal{P}_n) \subset \mathcal{P}_{n-\omega_{\ell'}}+\ldots+ \mathcal{P}_{n-1}+ \mathcal{P}_n.
\]
When computing $\nabla_\mathfrak{X}(\mathcal{P}_n)$,
the multiplication by $\zeta_{\ell,\ell'}\in\mathcal{P}_1$ is  then not able to compensate for the loss when $\omega_{\ell'}\geq 3$.
The profound reason is that our knowledge about the way~$\zeta_{\ell,\ell'}$ vanishes at the origin is too weak.

\begin{dfn}
A coordinate system adapted to the commutator structure of the vector fields $\mathfrak{X}$ near a regular H\"ormander point $\underline{x}$
of step $m$ is called \textsl{well-adapted} if
\begin{equation}
\forall \ell\in\{1,\ldots,n_1\}, \quad \forall \ell'\in\{1,\ldots,q\},\qquad
|\nabla_\mathfrak{X}^\gamma \zeta_{\ell,\ell'}| \leq C_{\gamma}\rho^{(\omega_{\ell'}-1-|\gamma|)_+}
\end{equation}
in a neighborhood of $\underline{x}$.
A family of vector fields that satisfies a regular H\"ormander condition is called \textsl{well-structured} if it admits a  well-adapted coordinate system.
\end{dfn}

One can check that in a well-adapted coordinate system, the gauge automatically satisfies~\eqref{FIELDSHIGHDERGAUGE}.
\begin{thm}[If-theorem for arbitrary step $m$] \label{ABSTRACTSTEPmFIELDS}
Let us consider a family of vector fields and $\underline{x}\in\Omega$ a regular H\"ormander point of step~$m$.
Then for any \textbf{well} adapted coordinate system, the gauge $\rho$ satisfies~\eqref{FIELDSHIGHDERGAUGE}
in the neighborhood of $\underline{x}$.
\end{thm}

\paragraph{Proof.}
The key is to adapt the definition of $\mathcal{P}_n$ to capture the enhanced knowledge that we gained about~$\zeta_{\ell,\ell'}$.
Let us define $\widetilde{\mathcal{P}}_n$ as the subset of $C^\infty(\Omega)$ that consists of homogeneous
polynomials with smooth coefficients of $x_\ell$, $\zeta_{\ell,\ell'}$ and of the derivatives of $\zeta_{\ell,\ell'}$ for which we have estimates, \ie
\[
\sum_{\alpha,\beta} x_1^{\alpha_1} \ldots x_q^{\alpha_q}
\zeta_{1,1}^{\beta_{1,1}} \zeta_{1,2}^{\beta_{1,2}} \ldots \zeta_{q,q}^{\beta_{q,q}}
\left(\prod_{\gamma}
(\nabla_{\mathfrak{X}}^{\gamma} \zeta_{1,1})^{\delta_{\gamma;1,1}} \ldots (\nabla_{\mathfrak{X}}^{\gamma} \zeta_{n_1,q})^{\delta_{\gamma;n_1,q}}
\right)
f_{\alpha,\beta,\delta}(x)
\]
where $f_{\alpha,\beta,\delta}\in C^\infty(\Omega)$, $\gamma$ denotes multi-indices of length $|\gamma|\geq1$ and
\[
\sum_{1\leq i\leq q} \alpha_i \omega_i
+ \sum_{j=1}^q\sum_{j'=1}^q \beta_{j,j'} (\omega_{j'}-1)
+ \sum_{|\gamma|\geq1} \sum_{j=1}^{n_1}\sum_{j'=1}^q \delta_{\gamma;j,j'} (\omega_{j'}-1-|\gamma|)_+
 = n.
\]
Note that only the factors for which $\omega_{j'}-1-|\gamma|>0$ are significant; the others can simply be tossed into $f_{\alpha,\beta,\delta}$.
For $n\leq 0$, one sets again $\widetilde{\mathcal{P}}_n = C^\infty(\Omega)$. We also introduce the linear span
\[
\widetilde{\mathcal{P}}_n^+ = \sum_{m\geq n} \widetilde{\mathcal{P}}_m.
\]
Using the Leibnitz formula, $\nabla_\mathfrak{X} (x_\ell) \in \widetilde{\mathcal{P}}_{\omega_\ell-1}$
and $\nabla_\mathfrak{X}(\widetilde{\mathcal{P}}_n) \subset \widetilde{\mathcal{P}}_{n-1}^+$. 
One also has
\[
\nabla_\mathfrak{X}\rho = \frac{\nabla_\mathfrak{X}(\rho^w)}{w\rho^{w-1}} \in \rho^{-(w-1)} \cdot \widetilde{\mathcal{P}}_{w-1}^+
\]
and recursively (note that $\rho^w \in \widetilde{\mathcal{P}}_w$ allows one to convert $\widetilde{\mathcal{P}}_0$ into $\rho^{-w}\cdot\widetilde{\mathcal{P}}_w$)
\[
\nabla_\mathfrak{X}^\gamma \rho \in \sum_{n\geq1} \frac{\mathcal{P}_n^+}{\rho^{n+|\gamma|-1}}
\]
from which \eqref{FIELDSHIGHDERGAUGE} follows immediately.
\cqfd\par

\bigskip
The previous ``abstract'' theorem does not presume on the existence of a well-adapted coordinate system.
However, when $m\leq 3$, it can actually be made to work.

\begin{thm}\label{STEP3FIELDS}
Any family of vector fields that satisfies a regular H\"ormander condition of step $m\leq 3$ 
admits at least one well-adapted coordinate system. It is therefore well-structured.
\end{thm}

\paragraph{Proof.}
For $m=1$ and $2$, any adapted coordinate system is well-adapted. Let us thus focus on~$m=3$ and use the previous notations.
Writing down the Taylor expansion of the coefficients for $\ell\leq n_1$:
\[
\zeta_{\ell,\ell'}(x) = \sum_{i\leq n_1} \left(\frac{ \partial \zeta_{\ell,\ell'} }{\partial x_i}(\underline{x})\right) x_i + O(\rho^2),
\]
it appears that, for $m=3$, a coordinate system is well-adapted if and only if
\begin{equation}\label{WELLADAPTEDm=3}
\forall \ell_1,\ell_2 \in\{1,\ldots,n_1\}, \quad \forall \ell_3\in\{n_2+1,\ldots,q\},\qquad
\frac{ \partial \zeta_{\ell_1,\ell_3} }{\partial x_{\ell_2}}(\underline{x}) = 0.
\end{equation}
Let us compute the following commutator:
\[
[Y_{\ell_1},Y_{\ell_2}] = \sum_{\ell=1}^q \left(
\frac{\partial \zeta_{\ell_2,\ell}}{\partial x_{\ell_1}} - \frac{\partial \zeta_{\ell_1,\ell}}{\partial x_{\ell_2}}
\right) \partial_\ell.
\]
At the point $\underline{x}$, the terms corresponding to $\ell> n_2$ must belong to $W_2(\underline{x})$ and thus vanish, therefore:
\begin{equation}\label{GOODNEWSm=3}
\forall \ell_1,\ell_2 \in\{1,\ldots,n_1\}, \quad \forall \ell_3\in\{n_2+1,\ldots,q\},\qquad
\frac{ \partial \zeta_{\ell_1,\ell_3} }{\partial x_{\ell_2}}(\underline{x}) = \frac{ \partial \zeta_{\ell_2,\ell_3} }{\partial x_{\ell_1}}(\underline{x}).
\end{equation}
One can now define a new coordinate system $(y_\ell)_{1\leq \ell\leq q}$ whose dual basis satisfies
\[
\frac{\partial}{\partial y_\ell} = \frac{\partial}{\partial x_\ell} +
\mathbf{1}_{\ell\leq n_1} \sum_{\ell'>n_2} \sum_{i\leq n_1} 
\left( \frac{ \partial \zeta_{\ell,\ell'} }{\partial x_i}(\underline{x})\right)  x_i  \cdot
\frac{\partial}{\partial x_{\ell'}} \cdotp
\]
This coordinate system is (locally) well defined because the fields $\frac{\partial}{\partial y_\ell} $ commute with each other
thanks to~\eqref{GOODNEWSm=3}.
By construction, this coordinate system satisfies~\eqref{WELLADAPTEDm=3} and is therefore a well-adapted one.
\cqfd\par

\begin{rmk}
The generalization of theorem~\ref{STEP3FIELDS} for $m\geq4$ is an \textsl{open} question. One can check that
a coordinate system is well-adapted if and only if
\begin{equation}
\frac{\partial^\alpha  \zeta_{\ell,\ell'}}{\partial x_1^{\alpha_1} \ldots \partial x_q^{\alpha_q}}(\underline{x}) = 0
\end{equation}
for any indices such that $\omega_\ell=1$, $\omega_{\ell'}\geq 3$ and $\sum_{i=1}^q \omega_i \alpha_i \leq \omega_{\ell'}-2$.
However, for $m\geq4$, it is not clear whether the regular H\"ormander assumption of step $m$ is enough to ensure
that the vector fields
\[
\frac{\partial}{\partial y_\ell} = \frac{\partial}{\partial x_\ell} +
\mathbf{1}_{\ell\leq n_1} \sum_{\ell'>n_2} \left(
\sum_{\sum \alpha_i\omega_i \leq \omega_{\ell'}-2}
\left( \frac{\partial^\alpha  \zeta_{\ell,\ell'}}{\partial x_1^{\alpha_1} \ldots \partial x_q^{\alpha_q}}(\underline{x}) \right)  x_1^{\alpha_1}\ldots x_q^{\alpha_q}  \cdot
\frac{\partial}{\partial x_{\ell'}}
\right)
\]
commute with each other.
\end{rmk}

\subsection{From the symbol property of the gauge to Hardy inequality}

For well-structured families of vector fields, symbols of class $S^k(\mathfrak{X};\rho)$ are functions $f$ such that
\begin{equation}
|\nabla_{\mathfrak{X}}^\gamma f(x) | \leq C_\gamma \rho(x)^{(k-|\gamma|)_+}
\end{equation}
in a neighborhood of $\underline{x}$, for any multi-index $\gamma$.
Once the symbol property is established for the gauge, the path that leads to the Hardy inequality is open.
The key (see \cite[chap.~7]{PhD}) is to define a ``radial'' vector field that admits both expressions:
\begin{equation}
R(x) = \sum_{\ell=1}^q \sigma_\ell(x) Y_\ell(x) = \sum_{k=1}^q (\omega_k x_k + \widetilde{\sigma}_k(x) )\partial_k
\end{equation}
in well-adapted coordinates, with $\sigma_\ell\in S^{\omega_\ell}(\mathfrak{X};\rho)$ and $\widetilde{\sigma}_k\in S^{\omega_k+1}(\mathfrak{X};\rho)$.
One can then check that
\begin{equation}
\div R = Q + O(\rho) \quad\text{and}\quad \lambda = \frac{R\rho}{\rho} \quad\text{satisfies}\quad
\begin{cases}
\lambda(x) = 1 + O(\rho),\\
R\lambda = O(\rho).
\end{cases}
\end{equation}
The computations of~\S\ref{PROOF} can then be carried out in a small enough neighborhood of $\underline{x}$.

\end{document}